\documentclass[reqno]{amsart}
\usepackage{amsmath}
\usepackage{amssymb}
\usepackage{color}

\newcommand\fuc{Fu\v c\'\i k\ }

\newtheorem{thm}{Theorem}[section]

\newtheorem{lemma}[thm]{Lemma}
\newtheorem{cor}[thm]{Corollary}
\newtheorem{prop}[thm]{Proposition}

\theoremstyle{definition}

\newtheorem{ex}[thm]{Example}

\theoremstyle{remark}
\newtheorem{remark}[thm]{Remark}
\newtheorem{remarks}[thm]{Remarks}

\numberwithin{equation}{section}

\newenvironment{mylist}
{\begin{enumerate}}
{\end{enumerate}}

\newcommand\al{\alpha}

\newcommand\ga{\gamma}
\newcommand\Ga{\Gamma}
\newcommand\de{\delta}
\newcommand\De{\Delta}

\newcommand\la{\lambda}
\newcommand\La{\Lambda}
\newcommand\om{\omega}

\newcommand\si{\sigma}
\newcommand\Si{\Sigma}
\renewcommand\th{\theta}

\newcommand\bal{{\boldsymbol{\alpha}}}
\newcommand\bfeta{\boldsymbol{\eta}}

\newcommand\br{{\boldsymbol{r}}}

\newcommand\bzero{{\boldsymbol{0}}}

\newcommand\tB{\widetilde B}
\newcommand\tX{\widetilde X}
\newcommand\tY{\widetilde Y}
\newcommand\tDe{\widetilde \Delta}

\newcommand\R{\mathbb{R}}

\newcommand\A{{\mathcal A}}
\newcommand\C{{\mathcal C}}

\newcommand\efun{u}

\newcommand\X{\times}

\newcommand{\pa}{\partial}
\newcommand{\sgn}{{\rm sgn\,}}
\renewcommand\emptyset{\text{\Large \o}}

\newcommand\beq{\begin{equation}}
\newcommand\eeq{\end{equation}}

\let\le=\leqslant
\let\ge= \geqslant

\begin{document}

\title[Multi-point problems]{Half eigenvalues and the Fu\v c\'\i k
spectrum of multi-point, boundary value problems}
\author{Fran\c cois Genoud, Bryan  P.  Rynne}
\address{Department of Mathematics and the Maxwell Institute for
Mathematical Sciences, Heriot-Watt University,
Edinburgh EH14 4AS, Scotland.}
\email{bryan@ma.hw.ac.uk}
\email{F.Genoud@ma.hw.ac.uk}

\begin{abstract}
We consider the nonlinear boundary value problem
consisting of the equation
\begin{equation}   \tag{1}
-u'' = f(u) + h,  \quad \text{a.e. on $(-1,1)$,}
\end{equation}
where $h \in L^1(-1,1)$,
together with the multi-point, Dirichlet-type boundary conditions
\begin{equation}   \tag{2}
u(\pm 1) = \sum^{m^\pm}_{i=1}\al^\pm_i u(\eta^\pm_i) ,
\end{equation}
where
$m^\pm \ge 1$ are integers,
$\al^\pm  = (\al_1^\pm, \dots,\al_m^\pm) \in [0,1)^{m^\pm}$,
$\eta^\pm \in (-1,1)^{m^\pm}$,
and we suppose that
$$
 \sum_{i=1}^{m^\pm} \al_i^\pm < 1 .
$$
We also suppose that $f : \R \to \R$ is continuous,
and
$$
0 < f_{\pm\infty} := \lim_{s \to \pm\infty} \frac{f(s)}{s} < \infty
$$
(we assume that these limits exist).
We allow
$f_{\infty} \ne f_{-\infty}$ --- such a nonlinearity $f$ is said to be
{\em jumping}.

Related to (1) is the equation
\begin{equation}   \tag{3}
-u''  = \la(a u^+ -  b u^-) ,
\quad \text{on $(-1,1)$,}
\end{equation}
where
$\la,\,a,\,b > 0$,
and $u^{\pm}(x) =\max\{\pm u(x),0\}$ for $x \in [-1,1]$.
The problem (2)-(3) is `positively-homogeneous' and jumping.
Regarding $a,\,b$ as fixed, values of $\la = \la(a,b)$ for which (2)-(3)
has a
non-trivial solution $u$ will be called {\em half-eigenvalues},
while the corresponding solutions $u$ will be called
{\em half-eigenfunctions}.

We show that a sequence of half-eigenvalues exists,
the corresponding half-eigenfunctions having specified nodal properties,
and we obtain certain spectral and degree theoretic properties
of the set of half-eigenvalues.
These properties lead to
solvability and non-solvability results for
the problem (1)-(2).
The set of half-eigenvalues
is closely related to the `\fuc spectrum' of the problem,
which we briefly describe.
Equivalent solvability and non-solvability results for (1)-(2)
are obtained from either the half-eigenvalue or the \fuc spectrum
approach.
\end{abstract}

\maketitle

\section{Introduction}  \label{intro.sec}

We consider the nonlinear boundary value problem
consisting of the equation
\begin{equation}  \label{orig.eq}
-u'' = f(u) + h,  \quad \text{a.e. on $(-1,1)$,}
\end{equation}
where $h \in L^1(-1,1)$,
together with the multi-point, Dirichlet-type boundary conditions
\begin{equation}  \label{dbc.eq}
u(\pm 1) = \sum^{m^\pm}_{i=1}\al^\pm_i u(\eta^\pm_i) ,
\end{equation}
where
$m^\pm \ge 1$ are integers,
$\al^\pm  = (\al_1^\pm, \dots,\al_m^\pm) \in \R^{m^\pm}$,
$\eta^\pm \in (-1,1)^{m^\pm}$.
We suppose that $f : \R \to \R$ is continuous,
and\beq  \label{f_limits.eq}
0 < f_{\pm\infty} := \lim_{s \to \pm\infty} \frac{f(s)}{s} < \infty
\eeq
(we assume that these limits exist).
We allow
$f_{\infty} \ne f_{-\infty}$ in  \eqref{f_limits.eq}
| such a nonlinearity $f$ is said to be {\em jumping}.

We will require some further notation, and conditions, for the
coefficients $\al^\pm$, which we now describe.
For any integer $m \ge 1$ and any point $\al \in \R^m$,
the notation $\al = 0$, $\al \ge 0$, $\al > 0$ will mean
$\al_i = 0$, $\al_i \ge 0$, $\al_i > 0$, for $i = 1,\dots,m$,
respectively;
$\A^m $ will denote the set of $\al \in \R^m $ satisfying
\begin{equation}  \label{A_cond.eq}
 \sum^{m}_{i=1} |\al_i| < 1 ,
\end{equation}
while
$\A_+^m $ will denote the set of $\al \in \A^m $ satisfying
$\al \ge 0$.
Throughout the paper we will suppose that $\al^\pm \in \A_+^{m^\pm}$.

Related to \eqref{orig.eq} is the equation
\begin{equation}   \label{fuc_de.eq}
-u''  = a u^+ -  b u^- ,
\quad \text{on $(-1,1)$,}
\end{equation}
where
$a,\,b > 0$,
and $u^{\pm}(x) =\max\{\pm u(x),0\}$ for $x \in [-1,1]$.
Putting $a = f_\infty$, $b = f_{-\infty}$, and $h=0$,
we can regard \eqref{fuc_de.eq} as a limiting form of \eqref{orig.eq} as
$|u| \to \infty$ (this will be made more precise below).
However, we will consider \eqref{fuc_de.eq}, and variants
of this equation, in its own right with general values of $a$ and $b$
(see Remark~\ref{abla_pos.rem}~(ii) below for why we assume that
$a,\,b > 0$).

The boundary value problem \eqref{dbc.eq}, \eqref{fuc_de.eq}
has a positively-homogeneous jumping nonlinearity,
in the sense that if $u$ is a solution then
$t u$ is also a solution for all $t \ge 0$.
We can now define  the {\em \fuc spectrum} of
\eqref{dbc.eq}, \eqref{fuc_de.eq}, to be the set
$$
\Si_F := \{ (a,b) \in \R^2 : \text{ \eqref{dbc.eq}, \eqref{fuc_de.eq}
has a non-trivial solution $u$}\}.
$$
The \fuc spectrum of problems with separated boundary conditions has
been used extensively to derive criteria for the solvability,
or non-solvability, of the general nonlinear problem
\eqref{orig.eq}, \eqref{dbc.eq}.
Such criteria have usually been described in terms of the location of
the point $(f_{\infty},f_{-\infty})$ in $\R^2$ relative to the spectrum
$\Si_F$.
See, for example, \cite{DAN1,DRA,RYN1}
and the references therein
(there are over 80 papers on Mathscinet with `\fuc spectrum' in the
title).

An alternative, real-valued, spectrum, which also yields
solvability criteria for
\eqref{orig.eq}, \eqref{dbc.eq},
can be defined by considering the  equation
\begin{equation}   \label{hev_de.eq}
-u''  = \la (a u^+ -  b u^-) ,
\quad \text{on $(-1,1)$,}
\end{equation}
with a spectral parameter $\la > 0$.
The problem \eqref{dbc.eq}, \eqref{hev_de.eq}
again has a positively-homogeneous jumping nonlinearity.
Regarding $a,\,b$ as fixed, we will say that any $\la$ for which
\eqref{dbc.eq}, \eqref{hev_de.eq}
has a non-trivial solution $u$ is a {\em half-eigenvalue},
and $u$ is a {\em half-eigenfunction},
and we define a corresponding spectrum to be the set
$$
\Si_H = \Si_H(a,b) :=
\{ \la \in \R : \text{ \eqref{dbc.eq}, \eqref{hev_de.eq}
has a non-trivial solution $u$}\}.
$$
For problems with separated boundary conditions the half-eigenvalue
spectrum has also been extensively investigated,
see for example
\cite{RYN1,RYN2},
and the references therein,
and a detailed comparison of the \fuc and half-eigenvalues approaches
to the solvability of
\eqref{orig.eq}, \eqref{dbc.eq}
is given in \cite{RYN2}.

Of course, in the current setting (with $a$, $b$ constant)
these spectra are closely related,
since a point $(a,b) \in \Si_F$ if and only if $1 \in \Si_H(a,b)$.
However, when the coefficients $a$, $b$ are variable
(which one would naturally consider when $f$, and the
limits $f_{\pm\infty}$, depend on $x$)
it is not so clear how to even define the \fuc spectrum,
and it is more difficult to obtain solvability criteria for the
general problem
\eqref{orig.eq}, \eqref{dbc.eq}.
In this case such criteria are usually given in terms of the location
of the set of points
$\{(f_\infty(x),f_{-\infty}(x)) : x \in (-1,1)\} \subset \R^2$,
relative to the constant coefficient \fuc spectrum $\Si_F$,
see \cite{RYN2} for  more details.
On the other hand, when $a$, $b$ are variable the spectrum
$\Si_H(a,b)$
can be defined just as above,
and it is shown in \cite{RYN2} that in this case the half-eigenvalue
approach can yield stronger solvability results than the \fuc spectrum
approach.

In our discussion of the multi-point problem here we will only
consider the constant coefficient problem, so the two approaches are,
in principle, equivalent.
Despite this, we will use the half-eigenvalue approach since our
method of obtaining the spectrum, and its properties,
will rely heavily on having a single parameter $\la \in \R$,
rather than a point $(a,b) \in \R^2$.

\begin{remark} \label{alt_form_hevals_de.rem}
The papers \cite{RYN1,RYN2} actually use half-eigenvalues
defined by equations of the form
\begin{equation} \label{alt_form_hevals_de.eq}
- u'' = a u^+ - b u^- + \la u ,
\quad \text{on $(-1,1)$,}
\end{equation}
rather than  \eqref{hev_de.eq}.
However, for the separated problem there is little difficulty
in converting the results of
\cite{RYN1,RYN2} to the setting of \eqref{hev_de.eq}.
Somewhat surprisingly, we find that this is not true in the
multi-point setting.
In fact, we find that it is considerably easier to deal
with the half-eigenvalue formulation \eqref{hev_de.eq} rather than
\eqref{alt_form_hevals_de.eq}, which is why we use \eqref{hev_de.eq}
here.
Indeed, the paper \cite{RYN4} considered \eqref{alt_form_hevals_de.eq}
together with the 3-point boundary condition
$$
u(-1) = 0, \quad u(1) = \al u(\eta)
$$
(in the present notation),
and encountered considerable technical difficulties.
The results we obtain here will generalise the results of
\cite{RYN4} to the much more general multi-point conditions
\eqref{dbc.eq}, while avoiding many of the difficulties
encountered in \cite{RYN4}.
\end{remark}

\begin{remark}
When $\al^\pm = 0$ the boundary conditions \eqref{dbc.eq} reduce
to the standard, separated Dirichlet conditions at $x = \pm 1$,
so we have termed the conditions \eqref{dbc.eq} `Dirichlet-type'.
Similar results to those below can also be obtained for multi-point
`Neumann-type' conditions, or a mixture of these,
see  \cite{GR} or \cite{RYN6}
(the Neumann-type problem is slightly more complicated because a
`multi-point operator' $\De$, which we introduce in
Section~\ref{f_space_De.sec} below, is non-singular since constant
functions lie in its kernel --- however, this can be dealt with in the
manner described in \cite{RYN6}).
More general nonlocal boundary conditions can also be considered, as
described in Section~1.1 of \cite{GR}.
For brevity we will not consider any of these cases further here ---
given the methods here, the extensions to these other cases is
relatively straightforward.
\end{remark}

We conclude this section with a brief description of the contents of the
paper.
Some preliminary results are described in Section~\ref{prelim.sec}.
In Section~\ref{heval.sec} we state (without proofs) our main
results on the spectral theoretic properties of the half-eigenvalue
problem, and on solvability and non-solvability conditions for a
corresponding inhomogeneous form of
\eqref{dbc.eq}, \eqref{hev_de.eq}
(the solvability conditions are expressed in terms of the
half-eigenvalues).
These results are similar to those obtained in \cite{RYN1} for separated
boundary conditions.
We also give some examples which show that our assumption that the
coefficients $\al^\pm \in \A_+^{m^\pm}$ is optimal,
and cannot be relaxed to the condition $\al^\pm \in \A^{m^\pm}$,
which is known to suffice for the linear, multi-point eigenvalue
problem (see \cite{RYN5}).
These results also enable us to easily construct the \fuc spectrum
$\Si_F$, and we give a brief description of this in
Section~\ref{fuc_spec.sec}.
We then begin the proofs.
In Section~\ref{single_bc.sec} we prove an existence and uniqueness
result for a problem consisting of equation \eqref{dbc.eq} together with
a single, multi-point, boundary condition
(this problem could be regarded as a multi-point analogue of the usual
initial value problem for equation \eqref{dbc.eq}).
As usual, the uniqueness result for this multi-point, `initial value
problem' then implies the simplicity of the half-eigenvalues
(in a suitable sense).
The results stated in Section~\ref{heval.sec} are then proved in
Sections~\ref{heval_proof.sec} and~\ref{hlslvble_proof.sec}.
In Section~\ref{nonlin.sec} we extend the
solvability and non-solvability results of
Section~\ref{heval.sec} to the more general problem
\eqref{orig.eq}, \eqref{dbc.eq},
and in Section~\ref{nonlin_prob.sec} we  obtain a global
bifurcation theorem, and nodal solutions of this problem.

\section{Preliminary definitions and results}  \label{prelim.sec}

\subsection{Function spaces and notation}  \label{f_space_De.sec}

For any integer $n \ge 0$, $C^n[-1,1]$ will denote the usual Banach
space
of $n$-times continuously differentiable functions on $[-1,1]$,
with the usual sup-type norm, denoted by $|\cdot|_n$
(all function spaces are taken to be real).
Let
$$
X := \{ u \in C^2[-1,1] : \text{$u$ satisfies \eqref{dbc.eq}} \} ,
\quad
Y := C^0[-1,1],
$$
with the norms $|\cdot|_2$ and $|\cdot|_0$ respectively.
We define a bounded linear operator $\De  : X \to Y$ by
$$
\De  u := -u'',  \quad  u \in X.
$$
This operator has a bounded inverse $\De ^{-1} : Y \to X$,
see \cite[Theorem~3.1]{RYN5}.
The following spaces will also be useful
$$
\tX := \{ u \in W^{2,1}(-1,1) : \text{$u$ satisfies \eqref{dbc.eq}} \} ,
\quad
\tY := L^1(-1,1) ,
$$
where $W^{2,1}(-1,1)$ denotes the usual Sobolev space of
functions with second order derivatives in $L^1(-1,1)$
($L^1(-1,1)$ and $W^{2,1}(-1,1)$ will be endowed with their usual
norms).
It can readily be checked that $\De $ extends to a bounded, invertible
operator $\tDe : \tX \to \tY$,
see \cite[Remark~3.4]{RYN5}.

Up to now we have regarded $\al^\pm,\,\eta^\pm,$ as constant,
and omitted them from the notation, for example, for the operator
$\De$.
However, in some of the discussion below it will be convenient to
regard some, or all, of these as variable,
and we will then indicate the dependence on these variables
in the obvious manner.
In particular, we let
$\bal := (\al^-,\al^+)$ and $\A_+ := \A_+^{m^-} \X \A_+^{m^+}$,
$\bfeta := (\eta^-,\eta^+)$,
and we will write, for example, $\De^{-1}(\bal)$ for $\bal \in \A_+$
(note that $X$ also depends on $\bal$, but $Y$ does not).

\subsection{Nodal properties}

We now introduce some notation to describe the nodal properties of
solutions of \eqref{orig.eq} or \eqref{fuc_de.eq}.
Firstly, for any $C^1$ function $u$, if $u(x_0)=0$ then $x_0$ is a
{\em simple} zero of $u$ if $u'(x_0) \ne 0$.
Now, for any integer $k \ge 1$ and any $\nu \in \{\pm\}$,
we define
$T_{k,\nu} \subset X$ to be the set of functions
$u \in X$ satisfying the following conditions:\\[1 ex]
T-(a) \ $u'(\pm 1) \ne 0$ and $\nu u'(-1) > 0$;\\
T-(b) \ $u'$ has only simple zeros in $(-1,1)$, and has exactly
$k$ such zeros;\\
T-(c) \  $u$ has a zero strictly between each consecutive zero of $u'$.
\\[1 ex]
We also define $T_k := T_k^+ \cup T_k^-$.
\medskip

\begin{remarks}
(i) If $u \in T_{k,\nu}$ then $u$
has exactly one zero between each consecutive zero of $u'$,
and all zeros of $u$ are simple.
Thus, $u$ has at least $k-1$ zeros in $(-1,1)$,
and at most $k$ zeros in $(-1,1]$.
\\(ii)\ \ The sets $T_{k,\nu}$ are open in $X$ and disjoint.
\\(iii) The sets $T_{k,\nu}$ were introduced in \cite{RYN3}.
Similar, but slightly different, sets have been used to describe
the nodal properties of solutions in the separated boundary condition
case,
but it is shown in \cite{RYN3} that the sets $T_{k,\nu}$ are more
suitable for the multi-point problem.
\\(iv)\ \ The symbols $\pm$ are used in two different contexts in this
paper:
(a) they refer to the end points $\pm 1$ at which the boundary
conditions hold (cf. \eqref{dbc.eq});
(b) they refer to the sign properties of the nodal solutions $u$,
in particular, through part~T-(a) of the above definition of the sets
$T_{k,\nu}$.
To attempt to keep the distinction between these usages clearer, we use
superscript $\pm$ to denote the boundary condition context, and
subscript  $\pm$
to denote the sign context
(except in the usage $u^\pm$ for the positive and negative parts of $u$,
as in equation \eqref{hev_de.eq} --- this usage seems too
well-established in the literature to change it here).
\end{remarks}

\section{Half-eigenvalues and solvability properties}  \label{heval.sec}

\subsection{Half-eigenvalues}

In this section we consider the half-eigenvalue problem
\eqref{dbc.eq}, \eqref{hev_de.eq},
which we rewrite in the form
\begin{equation} \label{heval.eq}
 - \De  u =  \la ( a u^+ - b u^- ) ,
\quad u \in X .
\end{equation}
General results on half-eigenvalues, and associated solvability and
degree theoretic results,  are described in \cite{BER} and \cite{RYN1}
for Sturm-Liouville problems with separated boundary conditions.
In this section we describe similar results for the multi-point problem.

\begin{thm}  \label{hevals.thm}
Suppose that $a,b > 0$ and $\bal \in \A_+$.
For each $k \ge 1$ and $\nu \in \{\pm\}$ there is exactly one
half-eigenvalue
$\la_{k,\nu} = \la_{k,\nu}(a,b) > 0$ of \eqref{heval.eq}
with a half-eigenfunction
$\efun_{k,\nu} = \efun_{k,\nu}(a,b) \in T_{k,\nu}$,
and there are no other half-eigenvalues.
The sequence of half-eigenvalues is strictly increasing, in the sense
that
\begin{equation}  \label{bermonot.eq}
k' > k \implies \la_{k',\nu'} > \la_{k,\nu},
\quad \text{for each $\nu',\,\nu \in \{\pm\}$}.
\end{equation}
In addition, $\pm \efun_{1,\pm}$ is strictly positive on $(-1,1)$,
and $\lim_{k \to \infty} \la_{k,\pm} = \infty$.
\end{thm}

Naturally, the half-eigenvalues in Theorem \ref{hevals.thm} depend on
all the parameters in the problem \eqref{heval.eq}.
The following result shows that this dependence is continuous.

\begin{cor} \label{cts_evals.cor}
The half-eigenvalues $\la_{k,\pm}$, $k \ge 1$, are $C^1$ functions of
the variables $(a,b,\bal,\bfeta)$
in $\R^2 \X \A_+ \X (-1,1)^{m^- + m^+} $.
\end{cor}

\begin{remarks} \label{abla_pos.rem}
It can easily be seen, from the form of the differential equation
\eqref{hev_de.eq} when $u$ is positive or negative, that if a
(non-trivial) solution $u$ changes sign at least 3 times then
$a,b,\la$
must all be non-zero and have the same sign,
so Theorem~\ref{hevals.thm} cannot hold unless this is true.
Thus, without loss of generality, we assume here that $a,b,\la$
are all positive.
\end{remarks}

\begin{remarks} \label{sep_hevals.rem}
For the separated half-eigenvalue problem
($a \ne b$, $\bal  = 0$),
similar results to those of Theorem~\ref{hevals.thm} are proved
in \cite[Theorem 2]{BER}, \cite[Theorem 4]{BRO}
and \cite[Theorem~5.1]{RYN1},
although these papers mainly consider half-eigenvalue problems of the
form
\begin{equation} \label{alt_form_hevals.eq}
-\De u = a u^+ - b u^- + \la u
\end{equation}
with variable coefficients $a,b$
(recall Remark~\ref{alt_form_hevals_de.rem}).
\end{remarks}

\begin{remark}[Eigenvalues] \label{evals.rem}
When $a=b$ (with $\bal  \ne \bzero$) the problem \eqref{heval.eq}
reduces to the linear multi-point eigenvalue problem
\begin{equation} \label{eval.eq}
 - \De u = \la a u, \quad u \in X .
\end{equation}
The spectral properties of \eqref{eval.eq}, with $a=1$, are obtained
in \cite[Theorem 3.1]{RYN5},
and it is clear that if we denote the eigenvalues and
eigenfunctions of this problem
by $\la_k$ and $\efun_k$, $ k \ge 1$,
then
$\la_{k,\pm}(a,a) = \la_k/a$ and we may suppose that
$\efun_{k,\pm}(a,a) = \pm \efun_k$.
Thus Theorem~\ref{hevals.thm} above extends the linear
eigenvalue results in \cite[Theorem 3.1]{RYN5}
to the half-eigenvalue problem
(however, \cite{RYN5} also deals with the $p$-Laplacian
eigenvalue problem, and the results in \cite{RYN5}
are valid for all $\al^\pm \in \A^{m^\pm}$).
\end{remark}

It is shown in \cite{RYN5} that for the eigenvalue problem
(with $a=b=1$)
Theorem~\ref{hevals.thm} holds for all $\bal \in \A$, but need not be
true if
$\bal \not\in \A$.
The following examples show that for the half-eigenvalue problem
Theorem~\ref{hevals.thm} need not be true if
$\bal \in \A \setminus \A_+$,
that is, if $\al^\pm$ satisfy \eqref{A_cond.eq} but some of the
coefficients $\al^\pm_i$ are negative.
Each of these examples have a Dirichlet condition at $-1$, so the
corresponding solutions can easily be sketched, and the truth of the
assertions should then be clear.

\begin{ex} \label{boundary.ex}
The problem
$$
\begin{array}{c}
-u'' = \left(\frac{2\pi}{3}\right)^2 u^+ - (\pi)^2 u^-
\\[1 ex]
u(-1) = 0, \quad u(1) = - \frac23 u(-\frac14) .
\end{array}
$$
has a half-eigenfunction $\efun \in \pa T_{1,+}$ given by
$$
\efun(x) :=
\begin{cases}
\frac{3}{2\pi} \sin \frac{2\pi}{3} (x+1) , \quad x \in (-1,\frac12) ,
\\[1 ex]
- \frac{1}{\pi} \cos \pi (x-1) , \quad x \in (\frac12,1) .
\end{cases}
$$
\end{ex}

\begin{ex} \label{missing.ex}
If $\de > 0$ is sufficiently small then the problem
$$
\begin{array}{c}
-u'' = \left(\frac{\pi}{2} + \de \right)^2 u^+ -
\left(\frac{1}{\de}\right)^2 u^-
\\[1 ex]
u(-1) = 0, \quad u(1) = - \frac12 u(0) .
\end{array}
$$
has no half-eigenfunction $\efun \in T_{2,+}$.
\end{ex}

\begin{remark}
It will be seen that the proofs of Theorems~\ref{hevals.thm}
and~\ref{branches.thm} below rely on the fact that eigenfunctions cannot
lie in the boundaries of the sets $T_{k,\nu}$
(so that the nodal properties of the eigenfunctions are preserved as
parameters are varied).
It will be shown that this is true when $\bal \in \A_+$, but
Example~\ref{boundary.ex} shows
that it need not be true when $\bal \in \A \setminus \A_+$.
\end{remark}

\subsection{Solvability properties}

In addition to eigenvalues, linear spectral theory is also concerned
with the solvability of inhomogeneous problems.
Accordingly, we will consider the
solvability of the inhomogeneous equation
\begin{equation}  \label{hlslvble.eq}
 - \tDe  u =  \la ( a u^+ - b u^- ) + h,
\quad  u \in \tX,
\end{equation}
for general functions $h \in \tY$, when $\la$ is not a
half-eigenvalue
(if $h \in Y$, then we would use the operator $\De$ and search for
solutions $u \in X$, in the same manner).

Clearly, \eqref{hlslvble.eq} is equivalent to the equation
\begin{equation}  \label{R_la_def.eq}
R_{\la} (u) := u + \la \tDe ^{-1}(a u^+ - b u^-) = - \tDe ^{-1} h,
\quad u \in \tX ,
\end{equation}
and the operator
$R_{\la} : \tX \to \tX$ is positively homogeneous, in the sense that
$R_{\la}(t u) = t R_{\la} (u)$ for any $t \ge 0$ and $u \in \tX$.
Also, the mapping
$u \to a u^+ - b u^- : \tX \to \tY$
is compact,
and hence the operator $I-R_{\la}$ is compact.
Now, letting $\tB_r(c)$ denote the ball in $\tX$ with
radius $r$ and centre $c$
(and putting $\tB_r := \tB_r(0)$, for brevity),
we see that the Leray-Schauder degree,
$\deg(R_{\la},\tB_1,0)$, is well-defined
whenever $\la$ is not a half-eigenvalue,
see \cite{DEI}.
We will now relate the solvability properties of \eqref{hlslvble.eq}
and the degree $\deg(R_{\la},\tB_1,0)$ to the location of $\la$
relative to the set of half-eigenvalues.
To state this precisely we introduce some further notation.

For each $k \ge 1$, let
$\la_{k,\max} = \max \{ \la_{k,+}, \la_{k,-} \}$,
$\la_{k,\min} = \min \{ \la_{k,+}, \la_{k,-} \}$,
and define the open intervals
\begin{align*}
\La_k^0 &=
\begin{cases}
(\la_{k,\min},\la_{k,\max}), & \text{if $ \la_{k,\min} <
\la_{k,\max}$},\\
\emptyset,                   & \text{if $ \la_{k,\min} = \la_{k,\max}$},
\end{cases}\\
\La_k^1 &= (\la_{k,\max},\la_{k+1,\min}),
\quad \La_0^1 = (-\infty,\la_{1,\min}).
\end{align*}
Intuitively, Theorem~\ref{hevals.thm} says that when $a \ne b$ the term
$a u^+ - b u^-$ in equation \eqref{heval.eq} `splits apart' the linear
eigenvalue $\la_k$ into
a pair of half-eigenvalues $\la_{k,\pm}(a,b)$, and the interval
$\La_k^0$ is the gap between these half-eigenvalues.
It is possible for the half-eigenvalues $\la_{k,\pm}$ to coincide,
so this gap may be empty.
On the other hand, the inequality \eqref{bermonot.eq} says that
half-eigenvalues with different values of $k$ cannot coincide,
so the interval $\La_k^1$ between half-eigenvalues corresponding to
$k$ and $k+1$ is non-empty.
Also, all these intervals are disjoint and their union comprises the
whole of $\R$, except for the half-eigenvalues.
Furthermore, by continuity, $\deg(R_{\la},\tB_1,0)$ is constant on
any of the intervals
$\La_k^0$, $\La_k^1$.

Clearly, the intervals $\La_k^i$, $i=0,1$, $k \ge 1$, depend on $a,\,b$,
and when it is necessary to indicate this dependence explicitly
we will write $\La_k^i(a,b)$.

\begin{thm}  \label{hlslvble.thm}
{\rm (A)}\  If $\la \in \La_k^1$, for some $k \ge 0$, then$:$
\begin{mylist}
\item
$\deg(R_{\la},\tB_1,0) = (-1)^k;$
\item
for any $h \in \tY$, equation \eqref{hlslvble.eq} has a solution
$u \in \tX$.
\end{mylist}
{\rm (B)}\
If $\la \in \La_k^0$, for some $k \ge 1$, then$:$
\begin{mylist}
\item
$\deg(R_{\la},\tB_1,0) = 0;$
\item
there exists $h \in \tY$ such that equation \eqref{hlslvble.eq}
has no solution$;$
\item
there exists $h_b \in \tY$ such that, for any $h \in \tB_1(h_b)$,
\eqref{hlslvble.eq} has at least two solutions.
\end{mylist}
\end{thm}

\begin{remark}
For the separated problem, similar results to those of
Theorem~\ref{hlslvble.thm} are proved in \cite[Theorem~1.4]{RUF}
(for constant coefficients)
and in \cite[Theorem~5.1]{RYN1}
(for variable coefficients).
The theorem shows that when $\la$ is not a
half-eigenvalue then a `nonlinear Fredholm alternative' holds for
\eqref{hlslvble.eq}, in the sense that either:
\\[2 pt]
(a)\ there exists a solution $u$ for all $h \in \tY$,
\\[2 pt]
or
\\[2 pt]
(b)\ there is no solution for some $h \in \tY$ and at least two
solutions for other $h \in \tY$.
\\[2 pt]
Such an interpretation was described in, for instance
\cite[Corollary~6.1]{RUF}.
\end{remark}

\begin{remark}
The above results can be regarded as a generalization, to the
half-linear, multi-point problem, of standard results from the linear
spectral theory of Sturm-Liouville problems with
separated boundary conditions.
When $a=b$ the problem reduces to the linear case, so Theorem
\ref{hlslvble.thm} also covers the linear, multi-point problem.
Of course, in the linear case $\la_{k,\min} = \la_{k,\max}$
for all $k \ge 1$,
so the intervals $\La_k^0$ are empty
and part (B) of Theorem \ref{hlslvble.thm} has no analogue.
Furthermore, the degree $\deg(R_{\la},\tB_1,0)$ changes by
2 as $\la$ crosses a linear eigenvalue, which can be
regarded, heuristically, as crossing two coincident half-eigenvalues,
each of which contributes a change of 1.
\end{remark}

The solvability and non-solvability results in
Theorem~\ref{hlslvble.thm}
will be extended to the general nonlinear
problem \eqref{orig.eq}, \eqref{dbc.eq},
in Theorem~\ref{nlslvble.thm} below

\subsection{The \fuc spectrum} \label{fuc_spec.sec}

Using the above results we can now construct the \fuc spectrum $\Si_F$.
We merely give a brief description here.

For $k \ge 1$, $\nu \in \{\pm\}$ and $\th \in (0,\pi/2)$,
let
$$
\hat\br(\th) := (\sin \th, \cos \th) ,
\quad
\la_{k,\nu}(\th) := \la_{k,\nu}(\sin \th, \cos \th),
$$
and define the $C^1$ curve
$$
\si_{F,k,\nu} :=
\{ \la_{k,\nu}(\th) \hat\br(\th) : \th \in (0,\pi/2)\} .
$$
It is clear from the preceding results that
$$
\Si_F = \bigcup_{k,\nu} \, \si_{F,k,\nu} .
$$
Some of the standard properties of the \fuc spectrum can easily be
deduced from this characterisation.
In particular, for each $k \ge 2$,
$(\la_k,\la_k) \in \si_{F,k,\pm}$
(that is, the curves $\si_{F,k,\pm}$ intersect the diagonal line
$\{(x,x) : x \in \R\} \subset \R^2$ at the point $(\la_k,\la_k)$),
and curves corresponding to different values of $k$ do not
intersect.
Also, it follows readily from the Sturm comparison theorem that
$$
\lim_{\th \to 0} \la_{k,\pm}(\th) =
\lim_{\th \to \pi/2} \la_{k,\pm}(\th) = \infty,
$$
so that the curves $\si_{F,k,\pm}$ have horizontal and vertical
asymptotes.
Furthermore, we can obtain analogues of the solvability properties
in parts (A) and (B) of Theorem~\ref{hlslvble.thm} in the gaps
(in the plane $\R^2$) between the pairs of curves $\si_{F,k,\pm}$,
or between the curves with consecutive values of $k$ ---
for more details of such solvability properties in the \fuc setting see,
for example, \cite{RYN2} or any of the other cited references dealing
with the \fuc spectrum.

Another standard property of the \fuc spectrum is that the curves
$\si_{F,k,\pm}$ are monotonically decreasing in $\R^2$.
This is not so easy to prove here, with the parametrisation of the
curves in terms of the angle $\th$.
Since the main interest of the \fuc spectrum is in obtaining
solvability criteria, and we have obtained such criteria in
Theorem~\ref{hlslvble.thm}, there seems little need to investigate the
geometrical properties of $\Si_F$ any further here.

The \fuc spectrum was first introduced by Fu\v c\'\i k in \cite{FUC},
and the paper \cite{DAN1} contains a comprehensive investigation
of this spectrum and its application to jumping nonlinearity problems.
Many of the results here have analogues in \cite{DAN1} and, in varying
degrees of generality, in many other papers since then (for separated
problems).

\section{Solutions with a single boundary condition}
\label{single_bc.sec}

In this section we will construct solutions of equation
\eqref{fuc_de.eq} satisfying a single, multi-point boundary condition.
For notational convenience in constructing such solutions,
from now on (in this section and later) we will write
$\la = s^2$, $a = \ga_+^2$, $b = \ga_-^2$
(given our hypothesis that $\la,a,b >0$ this is possible),
and we consider the problem
\begin{gather}
 -u'' =  s^2(\ga_+^2u^+ - \ga_-^2  u^-)  ,  \quad \text{on $\R$},
\label{single_bc_de.eq}
\\
u(\eta_0) = \sum^{m}_{i=1}\al_i u(\eta_i) ,
\label{single_bc.eq}
\end{gather}
for arbitrary $\al \in \A_+^m$, $\eta_0 \in \R$ and $\eta \in \R^m$.
We can regard
\eqref{single_bc_de.eq}, \eqref{single_bc.eq}
as a `multi-point, initial value problem',
and we will prove the following existence and `uniqueness' result for
this problem.

\begin{thm}  \label{single_bc.thm}
For fixed $s,\ga_\pm > 0$, $m \ge 1$, $\al \in \A_+^m$,
$\eta_0 \in \R$ and $\eta \in \R^m$,
there exist functions $\psi_\pm \in C^2(\R)$ such that
$\pm \psi_\pm'(\eta_0) > 0$
and
the set of solution of \eqref{single_bc_de.eq}, \eqref{single_bc.eq},
has the form
$$
\{ C_+ \psi_+ : C_+ \ge 0 \}  \cup \{ C_- \psi_- : C_- \ge 0 \} .
$$
\end{thm}

\begin{remark} \label{single_bc_linear.rem}
Theorem~\ref{single_bc.thm} shows that if $\ga_+ \ne \ga_-$ then,
in general, the solution set of
\eqref{single_bc_de.eq}, \eqref{single_bc.eq},
consists of two `half-rays' spanned by the functions $\psi_\pm$.
However, if $\ga_+ = \ga_-$
(so that the problem is linear),
then the solution set must be a linear subspace,
so we have $\psi_+ = \psi_-$ and the solution set has the form
$\{ C \psi_+ : C \in \R\}$.
\end{remark}

\begin{proof}

Define $\Psi \in C^2(\R)$ to be  the
solution of the initial value problem
\begin{equation}  \label{jump_ivp.eq}
\begin{array}{c}
-\Psi'' =  \ga_+^2\Psi^+ - \ga_-^2  \Psi^- ,
\\[1 ex]
\Psi(0)=0, \quad \Psi'(0) = 1.
\end{array}
\end{equation}
Clearly, $\Psi$ has only simple zeros and, on any interval
where $\pm \Psi > 0$, it satisfies the equation
$-\Psi'' = \ga_\pm^2 \Psi$,
so the graph of $\Psi$ consists of a succession of positive and
negative, sinusoidal `bumps'.
Thus,  $\Psi$ has the form
\begin{equation}  \label{psi_form.eq}
\Psi(x) =
\pm \frac{1}{\ga_\pm} \sin (\ga_\pm x - \tau_{\pm}(x)) ,
\quad \text{when $\pm \Psi(x) > 0$,}
\end{equation}
where $\tau_{\pm}(x)$ is defined by
$$
\ga_\pm^{-1}\tau_{\pm}(x) := \max \{ z \le x : \Psi(z) = 0 \},
$$
that is,  $\ga_\pm^{-1}\tau_{\pm}(x)$ is the zero of $\Psi$ immediately
below
$x$.
Also,  $\Psi$ is periodic, with period
$$
p_\Psi := \frac{\pi}{\ga_+} + \frac{\pi}{\ga_-}  .
$$

Next, for any $(s,\de) \in (0,\infty) \X \R$ we define
$w(s,\de) \in C^2(\R)$
by
$$
w(s,\de)(x) := \Psi(sx - \de) , \quad x \in \R .
$$
It can easily be verified that any solution of \eqref{single_bc_de.eq}
must have the form $C w(s,\de)$,
for some $C \ge 0$ and $(s,\de) \in (0,\infty) \X \R$,
and $w(s,\de)$ satisfies the boundary condition \eqref{single_bc.eq}
if and only if
\begin{equation}  \label{single_bc_zeros.eq}
\Ga(s,\de,\al) :=
w(s,\de)(\eta_0) - \sum_{i=1}^{m}  \al_i w(s,\de)(\eta_i) = 0 .
\end{equation}
Thus, it suffices to consider the set of
solutions of \eqref{single_bc_zeros.eq}.
Clearly, the function
$\Ga : (0,\infty) \X \R \X \A_+^m  \to \R$ is $C^2$,
and we will denote the partial derivatives of $\Ga$ with respect to
$s$ and $\de$ by $\Ga_s$ and $\Ga_\de$.

\begin{lemma}  \label{single_bc_simple.lem}
For any
$(s,\de,\al) \in (0,\infty) \X \R \X \A_+^m$,
$$
\Ga(s,\de,\al) = 0 \implies  \Ga_\de(s,\de,\al) \ne 0 .
$$
\end{lemma}

\begin{proof}
Suppose that,
for some $(s,\de,\al) \in (0,\infty) \X \R \X \A_+^m$,
\begin{equation}  \label{single_bc_Ga_zeros_simple.eq}
\Ga(s,\de,\al) = \Ga_\de(s,\de,\al) = 0 .
\end{equation}
Suppose also, without loss of generality, that
$$
w(s,\de)(\eta_0) \ge 0 , \quad
w(s,\de)(\eta_i) \ge  0 , \quad 1 \le i \le p , \quad
w(s,\de)(\eta_i) < 0 , \quad p < i \le m .
$$
Then, by \eqref{psi_form.eq}, we can rewrite
\eqref{single_bc_Ga_zeros_simple.eq} as
\begin{equation}  \label{S_C_z.eq}
 S_0   =
\sum_{i=1}^p \al_i S_i
-
\sum_{i=p+1}^m \al_i \frac{\ga_+}{\ga_-} S_i ,
\quad
C_0  = \sum_{i=1}^m \al_i C_i ,
\end{equation}
where
\begin{align*}
S_i  &=
\begin{cases}
\sin(\ga_+ (s\eta_i-\de) - \tau_{+}(s\eta_i-\de)),
\quad & 0 \le i \le p,
\\
\sin(\ga_- (s\eta_i-\de) - \tau_{-}(s\eta_i-\de)),
 & p < i \le m,
\end{cases}
\end{align*}
and the terms $C_i$, $0 \le i \le m$, are defined similarly, by
replacing
$\sin$ with $\cos$.
Now, by \eqref{S_C_z.eq},
and the fact that $S_i \ge 0$, $i = 0,\dots, m$,
\begin{align*}
1 & = S_0^2 + C_0^2
 \le
 S_0 \sum_{i=1}^p \al_i S_i  +  |C_0| \sum_{i=1}^m \al_i |C_i|
\le  \sum_{i=1}^m \al_i \big( S_0 S_i + |C_0| |C_i| \big)
\\ &
\le  \sum_{i=1}^m \al_i \big( S_0^2 +  C_0^2 \big)^{1/2}
\big( S_i^2 +  C_i^2 \big)^{1/2}
 <  1 ,
\end{align*}
which shows that \eqref{single_bc_Ga_zeros_simple.eq} cannot hold.
\end{proof}

\begin{lemma} \label{wd_nonzero.lem}
For any
$(s,\de,\al) \in (0,\infty) \X \R \X \A_+^m$,
$$
\Ga(s,\de,\al) = 0 \implies w(s,\de)'(\eta_0) \ne 0 .
$$
\end{lemma}

\begin{proof}
Suppose that,
for some $(s,\de,\al) \in (0,\infty) \X \R \X \A_+^m$,
\begin{equation}  \label{single_bc_Ga_zero_wd_nonzero.eq}
w(s,\de)(\eta_0) = \sum_{i=1}^{m}  \al_i w(s,\de)(\eta_i) > 0
\quad \text{and} \quad w(s,\de)'(\eta_0) = 0
\end{equation}
($w(s,\de)(\eta_0) \ne 0$ since $w(s,\de)$ is non-trivial,
and the case $w(s,\de)(\eta_0) < 0$ is similar).
Then, by its repeating, sinusoidal form, the function
$w(s,\de)$ has a global max at $x=\eta_0$,
that is,
$w(s,\de)(x) \le w(s,\de)(\eta_0)$, for all $x \in \R$,
so by \eqref{single_bc_Ga_zero_wd_nonzero.eq}
and the assumption that $\al \in \A_+^m$,
$$
w(s,\de)(\eta_0) = \sum_{i=1}^m \al_i w(s,\de)(\eta_i)
<
w(s,\de)(\eta_0) .
$$
This contradiction proves the lemma.
\end{proof}

For any $s > 0$ and $\al \in \A_+^m$,
it follows from the above definitions that the function
$\Ga(s,\cdot,\al)$ is $ p_\Psi $-periodic,
so there are multiple zeros of \eqref{single_bc_zeros.eq}
which do not yield distinct solutions of the problem
\eqref{single_bc_de.eq}, \eqref{single_bc.eq}.
To remove these additional zeros and to make the domain of
$ \de$ compact, from now on we will regard $\de$ as lying in
the circle
(which we denote by $S^1$)
obtained from the interval $[0,p_\Psi]$ by identifying the points $0$
and $p_\Psi$,
and we then regard the domain of
$\Ga$ as $(0,\infty) \X S^1 \X  \A_+$
(clearly, $\Ga$ is still $C^2$).

For any fixed $s > 0$ the function $\Ga(s,\cdot,0)$ has exactly
two zeros
$\de_\pm(s,0) \in S^1 $,
and these zeros are simple
and may be labelled so that $\pm w(s,\de_\pm(s,0))'(\eta_0) > 0$.
Hence, by a simple continuation argument,
using Lemma~\ref{single_bc_simple.lem} and the implicit function
theorem, we see that $\Ga(s,\cdot,\al)$ has exactly two zeros
$\de_\pm(s,\al ) \in S^1 $ for all $\al \in \A_+^m$,
and these zeros depend continuously on $\al$.
Also, by Lemma~\ref{wd_nonzero.lem},
$$
\pm w(s,\de_\pm(s,\al))'(\eta_0) = \pm w(s,\de_\pm(s,0))'(\eta_0) > 0,
\quad \al \in \A_+^m ,
$$
so setting  $\psi_\pm := w(s,\de_\pm(s,\al))$
completes the proof of Theorem~\ref{single_bc.thm}.
\end{proof}

\begin{remark}
The variable $s$ has been fixed throughout this section, but the above
notation and Lemmas~\ref{single_bc_simple.lem}
and~\ref{wd_nonzero.lem} will be required in
Section~\ref{heval_proof.sec},
where $s$ will vary.
\end{remark}

\section{Proof of Theorem \ref{hevals.thm} and
Corollary~\ref{cts_evals.cor}}
\label{heval_proof.sec}

Since we have assumed that $a,b,\la > 0$ we may now rewrite the
half-eigenvalue problem \eqref{heval.eq} in the form
\begin{equation}  \label{heval_rw.eq}
- \De  u =  s^2 (\ga_+^2 u^+ - \ga_-^2 u^- ) ,
\quad u \in X
\end{equation}
(that is, we have put $\la = s^2$, $a = \ga_+^2$, $b = \ga_-^2$),
and we may apply the constructions in Section~\ref{single_bc.sec}
to \eqref{heval_rw.eq}.
In particular, we continue to use the solution $w(s,\de)$ of the
differential equation \eqref{single_bc_de.eq} defined there.
Substituting $w(s,\de)$ into the boundary conditions \eqref{dbc.eq}
now shows that a number $\la = s^2$ is a half-eigenvalue of
\eqref{heval_rw.eq} if and only if the pair of equations
\begin{equation}  \label{pair_bc_zeros.eq}
\Ga^\pm(s,\de,\al^\pm) :=
w(s,\de)(\pm 1) - \sum_{i=1}^{m^\pm}  \al_i^\pm w(s,\de)(\eta_i^\pm)
= 0 ,
\end{equation}
is satisfied, for some $\de \in \R$,
and then $w(s,\de)$ is a corresponding half-eigenfunction.
Thus, we will prove the theorem by considering the set of solutions
of \eqref{pair_bc_zeros.eq}.
As for the function $\Ga$ in Section~\ref{single_bc.sec}, we will
regard the domains of the functions
$\Ga^\pm$ as $(0,\infty) \X S^1 \X  \A_+$.
Again, it is clear that these  functions are $C^2$ on this domain.

The following proposition now proves the existence and uniqueness
of the half-eigenvalues.

\begin{prop} \label{Zk_one_zero.prop}
For each $k \ge 1$, $\nu \in \{\pm\}$ and $\bal  \in \A_+ $,
there is exactly one solution
$(s_{k,\nu}(\bal),\de_{k,\nu}(\bal)) \in (0,\infty) \X S^1$
of \eqref{pair_bc_zeros.eq} such that
$w(s_{k,\nu}(\bal),\de_{k,\nu}(\bal)) \in T_{k,\nu}$.
There are no other solutions of \eqref{pair_bc_zeros.eq} in
$(0,\infty) \X S^1$.
\end{prop}

\begin{proof}
When $\bal = \bzero$ the problem \eqref{heval.eq} is a constant
coefficient, half-eigenvalue problem with separated (Dirichlet) boundary
conditions, so it is elementary to explicitly construct the
half-eigenvalues and corresponding half-eigenfunctions,
which we will write as
$\la_{k,\pm}^\bzero = (s_{k,\pm}^\bzero)^2$,
$\efun_{k,\pm}^\bzero = w(s_{k,\pm}^\bzero,\de_{k,\pm}^\bzero)
\in T_{k,\pm}$,
$k \ge 1$,
for suitable $s_{k,\pm}^\bzero$, $\de_{k,\pm}^\bzero$,
see \cite{DAN1} or the proof of \cite[Theorem 11.5]{DRA}
for the details.
This yields the following lemma (which we state for reference).

\begin{lemma} \label{bal_z.lem}
When $\bal = \bzero $, Proposition~\ref{Zk_one_zero.prop} holds
with corresponding solutions
$(s_{k,\nu}^\bzero,\de_{k,\nu}^\bzero)$,
$k \ge 1$, $\nu \in \{\pm\}$.
\end{lemma}

To extend  Lemma~\ref{bal_z.lem} to the case $\bal \ne \bzero$ we
require the following lemmas.

\begin{lemma} \label{Ga_nz_implies.lem}
For any
$(s,\de,\bal) \in (0,\infty) \X \R \X \A_+$, and $\nu \in \{\pm\}$,
$$
\Ga^\nu(s,\de,\al^\nu) = 0 \implies
\begin{cases}
w(s,\de)'(\nu) \ne 0,
\\[1 ex]
 \Ga_s^\nu(s,\de,\al^\nu) \, \Ga_\de^\nu(s,\de,\al^\nu) \ne 0
\end{cases}
$$
$($where $w(s,\de)'(\pm)$ means $w(s,\de)'(\pm 1))$.
\end{lemma}

\begin{proof}
This follows immediately from
Lemmas~\ref{single_bc_simple.lem} and~\ref{wd_nonzero.lem},
and  the definitions of $\Ga$ and $\Ga^\pm$
(the proof that $\Ga_s^\pm(s,\de,\al) \ne 0$ is similar to the proof
of Lemma~\ref{single_bc_simple.lem}, using the additional fact that
$|\eta_i^\pm| \le 1$, $i = 0,\dots, m^\pm$, here).
\end{proof}

\begin{lemma} \label{yla_in_Tk.lem}
For any
$(s,\de,\bal) \in (0,\infty) \X \R \X \A_+$,
$$
\Ga^+(s,\de,\al^+) = \Ga^-(s,\de,\al^-) = 0
\implies
\text{$w(s,\de) \in T_{k,\nu}$, for some $k \ge 1$ and $\nu \in
\{\pm\}$.}
$$
\end{lemma}

\begin{proof}
By the definition of the sets $T_{k,\nu}$ and the form of $w(s,\de)$,
it suffices to show that $w(s,\de)'(\pm 1) \ne 0$,
but this follows from Lemma~\ref{Ga_nz_implies.lem}.
\end{proof}

The proof of Proposition~\ref{Zk_one_zero.prop} can now be completed by
following the continuation argument in the proof of
\cite[Theorem~4.1]{GR},
so we will merely outline the argument here.
Lemmas~\ref{bal_z.lem}, \ref{Ga_nz_implies.lem}
and~\ref{yla_in_Tk.lem}
above provide the necessary analogues,
in the current setting,
of Corollary~4.5, Lemma~4.6 and equation (4.6) in \cite{GR}.
The argument is then similar to that in the proof of
Theorem~\ref{single_bc.thm}, except that we now have the pair of
equations \eqref{pair_bc_zeros.eq} to consider, rather than the single
equation \eqref{single_bc_zeros.eq}.
However, as shown in the the proof of
\cite[Theorem~4.1]{GR}, the above results enable us to apply the
implicit function theorem at an arbitrary solution $(s,\de,\bal)$
of \eqref{pair_bc_zeros.eq} (with $\bal = (\al^+,\al^-) \in \A_+$),
so we may construct the entire set of zeros of \eqref{pair_bc_zeros.eq},
at an arbitrary $\bal \in \A_+$,
by continuation away from the zeros at $\bal = \bzero$
(as given in Lemma~\ref{bal_z.lem}).
\end{proof}

Proposition~\ref{Zk_one_zero.prop} has proved the existence and
uniqueness of the half-eigenvalues
$\la_{k,\nu}(\bal) := s^2_{k,\nu}(\bal)$,
with half-eigenfunctions
$\efun_{k,\nu}(\bal) := w(s_{k,\nu}(\bal),\de_{k,\nu}(\bal))$
for $\bal \in \A_+$.
We will now prove that these half-eigenvalues are increasing,
in the sense of inequality \eqref{bermonot.eq}.
To do this we first note that, by the explicit construction when
$\bal = \bzero$, the half-eigenvalues
$\la_{k,\nu}(\bzero) = \la_{k,\nu}^\bzero$
satisfy \eqref{bermonot.eq}.
Thus, by the continuation construction of
$\la_{k,\nu}(\bal )$, for general  $\bal  \in \A_+ $,
it suffices to show that
$\la_{k,\nu}(\bal ) \ne \la_{k+1,\nu'}(\bal )$ for all
$k \ge 1$, $\nu,\nu' \in \{\pm\}$ and $\bal  \in \A_+ $.

Suppose, on the contrary, that
$\la_{k,\nu}(\bal) = \la_{k+1,\nu}(\bal )$,
for some such $k $, $\nu$ and $\bal $.
Then, by Proposition~\ref{Zk_one_zero.prop} and the definition of the
sets $T_{k,\nu},T_{k+1,\nu}$,
$$
\sgn \efun_{k,\nu}(\bal )'(-1) =  \sgn \efun_{k+1,\nu}(\bal )'(-1) ,
$$
so Theorem~\ref{single_bc.thm}, together with the boundary
condition \eqref{single_bc.eq} at $\eta_0 = -1$,
shows that
$\efun_{k,\nu}(\bal) = \efun_{k+1,\nu}(\bal)$,
which is a contradiction, since the sets $T_{k,\nu},T_{k+1,\nu}$
are disjoint.
Now suppose that
$\la_{k,\nu}(\bal ) = \la_{k+1,-\nu}(\bal )$.
Since, by definition, the derivative of a function in $T_{k,\nu}$
changes sign exactly $k$ times in the interval $(-1,1)$,
we see that in this case
$$
\sgn \efun_{k,\nu}(\bal )'(1) =  \sgn \efun_{k+1,-\nu}(\bal )'(1) ,
$$
so Theorem~\ref{single_bc.thm}, together with the boundary
condition \eqref{single_bc.eq} at $\eta_0 = 1$,
now shows that
$\efun_{k,\nu} = \efun_{k+1,-\nu}$,
which is again a contradiction.
This proves that the half-eigenvalues are increasing.

Next, the fact that $\lim_{k\to\infty} \la_{k,\pm} = \infty$ is clear
from the Sturm comparison theorem, so it only remain to prove that
$\pm\efun_{1,\pm}$ is strictly positive on $(-1,1)$.
To do this we observe that a half-eigenfunction which does not
change sign is in fact an eigenfunction of a linear problem.
Specifically, recalling Remark~\ref{evals.rem},
it is shown in \cite[Theorem 3.1]{RYN5} that the linear eigenfunction
$\efun_1$ may be chosen to be strictly positive on $(-1,1)$ and
$\efun_1 \in T_{1,+}$,
so by the above uniqueness result for the half-eigenvalues,
$$
\la_{1,+} = \la_1/a, \ \efun_{1,+} = \efun_1,
\quad
\la_{1,-} = \la_1/b, \ \efun_{1,-} = -\efun_1
$$
(up to a positive scalar multiple of the eigenfunctions).
This proves the desired positivity,
and finally completes the proof of Theorem~\ref{hevals.thm}.$\hfill\Box$

The proof of  Corollary~\ref{cts_evals.cor} now follows immediately from
the implicit function theorem and the continuation construction of
the half-eigenvalues,
using the fact that the functions $\Ga^\pm$ are $C^1$ functions of the
variables $(a,b,\bal,\bfeta)$
(although this dependence was suppressed above).

\section{Proof of Theorem \ref{hlslvble.thm}}
\label{hlslvble_proof.sec}

\subsection{Part (A)}

Corollary~\ref{cts_evals.cor} showed that the half-eigenvalues depend
continuously on $\bal  \in \A_+ $, and for the duration of this proof
we indicate this dependence explicitly by writing
$\la_{k,\pm}(\bal )$, $\bal  \in \A_+$, for $k \ge 1$.
Furthermore, for the duration of this proof, we indicate the dependence
on $\bal$ of the space $\tX$,
and hence the ball $\tB_1 \subset \tX$ and the operator
$R_\la : \tX \to \tX$,
by writing
$\tX_\bal $, $\tB_{1,\bal }$ and $R_{\la,\bal}$.
We can readily construct a continuous family of
bounded, linear isomorphisms
$S_\bal  : \tX_\bzero \to \tX_\bal $, $\bal  \in \A_+ $,
with $S_\bzero$ the identity on $\tX_\bzero$.
We now fix $k \ge 1$ and prove the result for arbitrary (fixed)
$ \bal^0 \in \A_+$ and
$\la^0 \in (\la_{k,\max}(\bal^0 ),\la_{k+1,\min}(\bal^0))$.

By \eqref{bermonot.eq}, we can choose a
continuous function
$\rho : [0,1] \to \R$ such that
\beq  \label{rho_in_La.eq}
 \rho(t) \in (\la_{k,\max}(t \bal^0 ),\la_{k+1,\min}(t \bal^0 )),
\quad t  \in [0,1] ,
\qquad \rho(1) = \la^0,
\eeq
and define $T_t  : \tX_\bzero \to \tX_\bzero$, $t  \in [0,1]$, by
$$
T_t(u) := S_{t\bal^0}^{-1} R_{\rho(t),t\bal^0} S_{t\bal^0}  u ,
\quad u \in \tX_\bzero   .
$$
Clearly, $T_t(u)$ depends continuously on
$(t ,u) \in [0,1] \X \tX_\bzero$.
Also, by \eqref{rho_in_La.eq} and the definition of $T_t$,
for each $t \in [0,1]$ there is no non-trivial solution
of the equation $T_t(u) = 0$,
so by standard properties of the degree,
$$
\deg(R_{\la^0},\tB_{1,\bal^0},0) =  \deg(T_1,\tB_{1,\bzero},0)
= \deg(T_0,\tB_{1,\bzero},0)
= \deg(R_{\rho(0),\bzero},\tB_{1,\bzero},0) .
$$
Now,
recalling that when $\bal = \bzero$ the problem is a Dirichlet
problem,
it is shown in \cite[Theorem~5.5]{RYN1}
(which deals with the separated problem)
that
$$
\deg(R_{\rho(0),\bzero},\tB_{1,\bzero},0) =  (-1)^k ,
$$
which proves part (A)-(a).
Part (A)-(b) now follows from part (A)-(a), the positive homogeneity
of the operator $R_{\la^0}$, and standard properties
of the degree, see parts (D4) and (D5) of Theorem 8.2 in \cite{DEI}.
Finally, a simple modification of the above argument also proves
the result for the case $k=0$
(simply by omitting any reference to $\la_{0,\max}$ in this case).

\subsection{Part (B)}

We first prove part (B)-(b) --- part (B)-(a) then follows immediately
since, by the above argument proving part A-(b),
the existence of $h$
for which equation \eqref{hlslvble.eq} has no solution shows that
$\deg(R_{\la},\tB_1,0) = 0$.

For any $s>0$,  we set $m = m^-$, $\eta_0 = -1$, $\eta = \eta^-$ and
$\al = \al^-$,
and use the results of Section~\ref{single_bc.sec}
to construct the corresponding numbers $\de_\pm(s,\al^-)$,
and  functions $\psi_\pm = w(s,\de_\pm(s,\al^-))$
which satisfy  \eqref{fuc_de.eq} and the boundary condition
\eqref{dbc.eq} at $-1$
(as in Theorem~\ref{single_bc.thm}).
Thus, recalling the functions  $\Ga^\pm$
defined in \eqref{pair_bc_zeros.eq},
in the proof of Theorem~\ref{hevals.thm},
we see that
$\Ga^-(s,\de_\pm(s,\al^-),\al^-) = 0$,
and $\la = s^2$ is a half-eigenvalue iff
$$
B(\la) :=
\Ga^+(s,\de_-(s,\al^-),\al^+)
\Ga^+(s,\de_+(s,\al^-),\al^+)
= 0 .
$$
Also, the results in the proof of Theorem~\ref{hevals.thm} show that if,
for some $k \ge 1$, $\la_{k,-} \ne \la_{k,+}$
then the sign of $B(\la)$ changes as $\la$ crosses the half-eigenvalues
$\la_{k,\pm}$.

\begin{prop} \label{Bla_gez.prop}
If  $B(\la) > 0$ then there exists $h \in \tY$ such that equation
\eqref{hlslvble.eq} has no solution.
\end{prop}

\begin{proof}
Suppose further that
\begin{equation}  \label{Ga_pos.eq}
\Ga^+(s,\de_\pm(s,\al^-),\al^+) > 0
\end{equation}
(the other case is similar).
In particular, \eqref{Ga_pos.eq} implies that
$\psi_\pm(1) > 0$.
We first construct a suitable function $h$.
Consider the initial value problem
\begin{equation} \label{subsidIVP.eq}
\begin{array}{c}
-v'' = a v + \lambda v - 1,
\\[1 ex]
v(x_l) = 0, \quad v'(x_l) \ge 0,
\end{array}
\end{equation}
for arbitrary $x_l \in (-1,1)$.
It can easily be shown that there exists a sufficiently small $\delta>0$
such that if $x_l \in [1-\delta,1)$
then any solution $v$ of \eqref{subsidIVP.eq}
has no zero in $(x_l,1]$.
We can also choose $\delta$ sufficiently small that
$\eta_i^\nu < 1-\delta$,
for $1 \le i \le m^\nu$, $\nu \in \{\pm\}$,
and
$\psi_\pm(x) \ne 0$
for $x \in [1-\delta,1]$.
Let $x_0 = 1-\de$, and define
$$
h(x) =
\begin{cases}
0,  &  x \in [0,x_0),
\\
-1,  &  x \in [x_0,1].
\end{cases}
$$
Now, for any $\ga \in \R$, let
$\Phi_{\ga,0} := |\ga| \psi_{\sgn \ga}$
(so $\Phi_{\ga,0}$ satisfies \eqref{fuc_de.eq}  and the boundary
condition \eqref{dbc.eq} at $-1$)
and let $\Phi_{\ga,h}$ denote the solution of
the differential equation corresponding to \eqref{hlslvble.eq}
satisfying $\Phi_{\ga,h} \equiv \Phi_{\ga,0}$ on $[-1,x_0]$.
Any solution of \eqref{hlslvble.eq} must be of the form
$\Phi_{\ga,h}$, for some $\ga \in \R$,
and $\Phi_{\ga,h}$ is a solution of \eqref{hlslvble.eq} if and only if
it satisfies the boundary condition \eqref{dbc.eq} at $1$.
We will show that this cannot happen for any $\ga \in \R$.

Clearly, if $\ga = 0$ then $\Phi_{\ga,0} \equiv 0$,
while if $\ga \ne 0$ then $\Phi_{\ga,0}$ has no zeros on $[x_0,1]$.
Furthermore, from the form of $h$,
$$
\Phi_{\ga,0}(x_0) = \Phi_{\ga,h}(x_0),
\quad
\Phi_{\ga,0}'(x_0) = \Phi_{\ga,h}'(x_0),
$$
and $\Phi_{\ga,h}(x) > \Phi_{\ga,0}(x)$,
for sufficiently small $x-x_0 > 0$.
In fact, we have the following result.

\begin{lemma}  \label{Phih_Phiz.lem}
For any $\ga \in \R$, $\Phi_{\ga,h} > \Phi_{\ga,0}$,
on $(x_0,1]$.
\end{lemma}

\begin{proof}
Suppose the contrary, and let $x_1$ be the first zero  of
$\Phi_{\ga,h} - \Phi_{\ga,0}$ in $(x_0,1]$.
Then $\Phi_{\ga,h} - \Phi_{\ga,0} > 0$ in $(x_0,x_1)$,
and $\Phi_{\ga,h}'(x_1) - \Phi_{\ga,0}'(x_1) \le 0$.
Now suppose that $\Phi_{\ga,0} \ge 0$ on
$(x_0,1)$, and define
$$
W =  \Phi_{\ga,h}' \Phi_{\ga,0} - \Phi_{\ga,h}  \Phi_{\ga,0}'.
$$
Then $W(x_0)=0$, and by equation \eqref{hlslvble.eq}, $W' > 0$ on
$(x_0,x_1)$,
so that $W > 0$ on $(x_0,x_1]$.
On the other hand,
$$
W(x_1) = \Phi_{\ga,0}(x_1) \bigl( \Phi_{\ga,h}'(x_1) -
\Phi_{\ga,0}'(x_1) \bigr) \le 0,
$$
and this contradiction deals with the case $\Phi_{\ga,0} \ge 0$ on
$(x_0,1)$.
Now suppose that $\Phi_{\ga,0} < 0$ on $(x_0,1)$.
If $\Phi_{\ga,h} < 0$ on $(x_0,x_1)$ then a similar argument deals
with this case so we suppose that $\Phi_{\ga,h}$ changes sign on
$(x_0,x_1)$.
But, by the choice of $\de$, the function $\Phi_{\ga,h}$ can
have at most one zero in $(x_0,1]$,
so we again obtain $\Phi_{\ga,h} > \Phi_{\ga,0}$,
on $(x_0,1]$,
which completes the proof of Lemma~\ref{Phih_Phiz.lem}.
\end{proof}

Now, for any $\ga \in \R$, combining \eqref{Ga_pos.eq} and
Lemma~\ref{Phih_Phiz.lem} yields
\begin{align*}
\Phi_{\ga,h}(1) &> \Phi_{\ga,0}(1)
= |\ga| \psi_{\sgn \ga}(1)
\ge  \sum_{i=1}^{m^+} \al_i^+ |\ga| \psi_{\sgn \ga}(\eta_i^+)
\\ &
=  \sum_{i=1}^{m^+} \al_i^+  \Phi_{\ga,0}(\eta_i^+)
=  \sum_{i=1}^{m^+} \al_i^+  \Phi_{\ga,h}(\eta_i^+) ,
\end{align*}
which completes the proof of Proposition~\ref{Bla_gez.prop}.
\end{proof}

Now, combining Proposition~\ref{Bla_gez.prop} with part (A) shows that
$$
\bigcup_{k \ge 1} \La_k^0 = \{ \la : B(\la) > 0 \},
\quad
\bigcup_{k \ge 0} \La_k^1 = \{ \la : B(\la) < 0 \},
$$
and hence completes the proof of part (B)-(b) of
Theorem~\ref{hlslvble.thm}.

To prove part (B)-(c) we require the following two lemmas.
The first lemma shows that the operator
$R_\la : \tX \to \tX$
has a Fr\'echet derivative at certain points $u \in \tX$,
which we denote by $D_u R_{\la}(u)$,
and gives the form of this derivative.
For any $u \in \tX$ we  denote the characteristic functions of the
sets $\{x \in [-1,1] : \pm u(x) > 0\}$ by $\chi_{u^\pm}$.

\begin{lemma}  \label{R_la_deriv.lem}
If $u \in \tX$ has only simple zeros in $[-1,1]$ then,
for any $\la > 0$, the operator $R_\la$ is
Fr\'echet differentiable at $u$, with $D_u R_{\la}(u)$  given by
$$
D_u R_{\la}(u) v
= v + \la \tDe^{-1} \big((a \chi_{u^+} + b \chi_{u^-})v\big),
\quad v \in \tX ,
$$
\end{lemma}

\begin{proof}
The proof is similar to the proof of  \cite[Lemma~3.1]{RUF}.
\end{proof}

\begin{lemma}  \label{om_non_sing.lem}
If $\la \in \La_k^0$, for some $k \ge 1$, then there exists
$ \om = \om(\la) \in \tX$
such that $D_u R_{\la}(\om)$ is non-singular.
\end{lemma}

\begin{proof}
Choosing $\om_1 \in \tX$ such that $\om_1 > 0$,
it follows from Lemma~\ref{R_la_deriv.lem} that
$$
D_u R_{\la}(\om_1) v = v + \la \tDe^{-1} (a v),
\quad v \in \tX ,
$$
and hence $D_u R_{\la}(\om_1)$ is non-singular
(so $\om_1$ suffices for the result)
unless the equation
\begin{equation} \label{non_sing_1.eq}
- \tDe v =  \la a v, \quad v \in \tX ,
\end{equation}
has a non-trivial solution $v$.

Suppose that \eqref{non_sing_1.eq} has a non-trivial solution
$v_1$.
Choose $x_j$, $j = 1,2$, such that $\eta_i^\nu < x_1 < x_2 \le 1$
for all $i = 1,\dots,m^\nu$, $\nu \in \{\pm\}$,
and choose  $\om_2 \in \tX$,
such that $\om_2$ has a simple zero at each $x_j$, and
$\om_2 > 0$ on $[-1,x_1) \cup (x_2,1]$,
$\om_2 < 0$ on $(x_1,x_2)$.
We now have
$$
D_u R_{\la}(\om_2) v =
  v + \la \tDe^{-1}(a \chi_{\om_2^+} + b \chi_{\om_2^-})v,
\quad v \in \tX
$$
(note that since $\la \in \La_k^0$ we must have $a \ne b$),
and now $D_u R_{\la}(\om_2)$ is non-singular
(so $\om_2$ suffices for the result)
unless the equation
\begin{equation} \label{non_sing_2.eq}
-\tDe v = \la (a \chi_{\om_2^+} + b \chi_{\om_2^-}) v, \quad v\in\tX ,
\end{equation}
has a non-trivial solution $v$.

Suppose that \eqref{non_sing_2.eq}  has a non-trivial solution $v_2$.
Given the choice of $x_1,\,x_2$, it follows from
the boundary condition \eqref{dbc.eq} at $x = -1$ and
Remark~\ref{single_bc_linear.rem},
together with
equations \eqref{non_sing_1.eq} and \eqref{non_sing_2.eq},
that we may suppose that
$ v_1 \equiv v_2 , \ \text{on $[-1,x_1]$,} $
and hence, by the boundary condition \eqref{dbc.eq} at $x=1$,
we must have
\begin{equation} \label{v1_eq_v2.eq}
v_1(1) = v_2(1) .
\end{equation}
We now show that there exist $x_1,x_2$ such that \eqref{v1_eq_v2.eq} is
false.

\begin{lemma} \label{sturm_comp.lem}
If $x_2 = 1$ and $x_1$ is sufficiently close to 1 then
$v_1(1) \ne v_2(1)$.
\end{lemma}

\begin{proof}
Suppose that $v_1(1) = v_2(1)$.
Then we may choose $x_1$ sufficiently close to 1 that $v_1,\,v_2$
are non-zero on $[x_1,1)$, and satisfy
\begin{gather*}
-v_1'' = \la a v_1 , \quad -v_2'' = \la b v_2
\\
v_1(x_1) = v_2(x_1) , \quad v_1'(x_1) = v_2'(x_1) , \quad
v_1(1) = v_2(1) .
\end{gather*}
Since $a \ne b$, a slight modification of the proof of Theorem 1.2 in
Chapter 8 of \cite{CL}
(the Sturm comparison theorem)
now shows that this is impossible.
\end{proof}

Returning to the proof of Lemma~\ref{om_non_sing.lem},
we choose $x_1$ as in Lemma~\ref{sturm_comp.lem},
and we can then choose $x_2 < 1$ so that $v_1(1) \ne v_2(1)$
(by continuous dependence of $v_2$ on $x_2$).
This shows that, for this choice of $x_1,x_2$, \eqref{v1_eq_v2.eq} is
false, and hence $D_u R_{\la}(\om_2)$ is non-singular.
This completes the proof of Lemma~\ref{om_non_sing.lem}.
\end{proof}

We now proceed with the proof of part (B)-(c) of
Theorem~\ref{hlslvble.thm}.
Let $h_b := - \tDe R_{\la}(\om) \in \tY$, where $\om$ is as in
Lemma~\ref{om_non_sing.lem}.
Then $\om$ is an isolated solution of the equation
$R_{\la}(u) = - \tDe^{-1} h_b$,
with index $\pm 1$.
Hence, by continuity properties of the degree and part (B)-(a), there
exists sufficiently small numbers $r_1,\,r_2,\,r_3 > 0$ such that,
for any $h \in \tB_{r_3}(h_b)$,
$$
\deg(R_{\la},\tB_1(0), - r_1 \tDe^{-1} (h) ) = 0,
\quad
\deg(R_{\la},\tB_{r_2}(r_1 \om), - r_1 \tDe^{-1} (h) )
\ne
0,
$$
and so the equation
$R_{\la}(u) = - r_1 \tDe^{-1} (h)$
has solutions in the balls $\tB_{r_2}(r_1\om)$ and
$\tB_1(0) \setminus \tB_{r_2}(r_1\om)$
(we assume that the numbers $r_i$ are sufficiently small that
$\tB_{r_2}(r_1\om) \subset \tB_1(0)$).
The result as stated in the theorem now follows by scaling and using the
positive homogeneity of $R_{\la}$.
This proves part (B)-(c), and so completes the proof of
Theorem~\ref{hlslvble.thm}.

\section{Non-linear problems}  \label{nonlin.sec}

In this section we consider the solvability properties
of the problem \eqref{orig.eq}, \eqref{dbc.eq},
which we rewrite as
\begin{equation}  \label{orig_rew.eq}
- \De  u = f(u) + h,
\quad  u \in \tX
\end{equation}
(here, for any $u \in C^0[-1,1]$ we let $f(u) \in C^0[-1,1]$ denote the
function $f(u(x))$, $x \in [-1,1]$).
We have the following analogue of
Theorem~\ref{hlslvble.thm}.

\begin{thm}  \label{nlslvble.thm}
{\rm (A)}\
If $1 \in \La_k^1(f_{\infty}, f_{-\infty})$, for some $k \ge 0$, then
for any $h \in \tY$
equation \eqref{orig_rew.eq} has a solution $u \in \tX$.
\\
{\rm (B)}\
If $1 \in \La_k^0(f_{\infty}, f_{-\infty})$, for some $k \ge 1$, then
there exists
$h_0,\,h_2 \in \tY$ such that if $h=h_0$ then
equation \eqref{orig_rew.eq} has no solution,
while if $h=h_2$ then
equation \eqref{orig_rew.eq} has two solutions.
\end{thm}

\begin{proof}
The proof of part (A) and the zero solution assertion in part (B)
is similar to the proof of parts
(iii) and (iv) of \cite[Theorem 5]{DAN1}, using the results of
Theorem~\ref{hlslvble.thm} above,
while the proof of the two solutions assertion is a slight extension of
the above proof of the corresponding result in
Theorem~\ref{hlslvble.thm}.
\end{proof}

The solvability results stated in
Theorem~\ref{nlslvble.thm}
were obtained in \cite[Theorem~5]{DAN1} for the equation
\eqref{orig_rew.eq} with Dirichlet boundary conditions.
These results were extended to more general, separated boundary
conditions and variable coefficients in \cite[Theorem~6.1]{RYN1}.
Similar results for other separated problems have
been obtained in many other papers
(see, for example,  the references in \cite{RYN1,RYN2}).
In many of these papers the results have been expressed in terms of
hypotheses on the location
of the point $(f_{\infty}, f_{-\infty}) \in \R^2$ relative to the
\fuc spectrum of the problem.
The hypotheses in Theorem~\ref{nlslvble.thm} ensure that
$\la=1$ is not a half-eigenvalue, and so $(f_{\infty}, f_{-\infty})$
is not in the \fuc spectrum.
In fact, in the constant coefficient case considered here
these hypotheses are equivalent to the usual conditions on the
location of the point $(f_{\infty}, f_{-\infty})$ relative to
the \fuc spectrum
(see \cite{RYN1,RYN2} for a more detailed discussion of the relationship
between the half-eigenvalues and the \fuc spectrum in the separated
problem, which applies equally well here).
That is, the \fuc spectrum and the half-eigenvalues are
equivalent concepts here.
However, it was convenient to prove Theorem~\ref{hevals.thm},
in particular, using the half-eigenvalue parameter $\la$.
Furthermore, as mentioned in the introduction,
it is shown in \cite{RYN1,RYN2} that in the general,
separated, variable coefficient case the results obtained using
half-eigenvalues are stronger than those obtained using the
\fuc spectrum approach, so it seemed preferable to state our results in
terms of half-eigenvalues rather than the \fuc spectrum.

\section{Global bifurcation and nodal solutions}
\label{nonlin_prob.sec}

In this final section we suppose that
\begin{equation}  \label{fzero.eq}
f(0) = 0, \quad f_0 := \lim_{s \to 0} \frac{f(s)}{s} > 0
\end{equation}
(we assume that this limit exists and is finite),
and we briefly describe a Rabinowitz-type global bifurcation theorem
and then obtain nodal solutions of \eqref{orig_rew.eq} with $h=0$.
Given the preceding results, these results are now relatively standard
so we omit most of the details.

\subsection{Global bifurcation}

We first briefly consider the bifurcation problem,
\begin{equation}  \label{bif.eq}
 - \De(u) =\la f(u),  \quad (\la,u) \in \R \X X .
\end{equation}
Clearly, by \eqref{fzero.eq}, $u \equiv 0$ is a solution of
\eqref{bif.eq} for any
$\la \in \R$;
such solutions will be called {\em trivial}.
The following Rabinowitz-type global bifurcation
theorem for nontrivial solutions of \eqref{bif.eq}
was proved in \cite[Theorem 6.2]{RYN5}
(recall the eigenvalues $\la_k$ introduced in Remark~\ref{evals.rem}).

\begin{thm}  \label{branches.thm}
For each $k \ge 1$ there exist closed, connected sets
$\C_k^\pm \subset (0,\infty) \X X$ of solutions of
\eqref{bif.eq} with the properties$:$
\begin{mylist}
\item
$(\la_k / f_0,0)\in {\mathcal C_k^\pm};$
\item
${\mathcal C_k^\pm}\backslash \{(\la_k / f_0,0)\} \subset
(0,\infty) \times T_k^\pm;$
\item
${\mathcal C_k^\pm}$ is unbounded in $(0,\infty) \X Y$.
\end{mylist}
\end{thm}

\subsection{Nodal solutions}

We now search for nodal solutions of the problem
\begin{equation}  \label{nodal.eq}
- \De  u = f(u) ,  \quad  u \in X
\end{equation}
(that is, solutions $u$ lying in specific sets $T_{k,\nu}$).
We suppose through this section that
\begin{equation}  \label{sfs_pos.eq}
s f(s) > 0, \quad s \ne 0 .
\end{equation}

\begin{thm} \label{jumping.thm}
If \eqref{sfs_pos.eq} holds and
\begin{equation}  \label{jumping_cond.eq}
(\la_k/f_0 - 1) ( \la_{k,\nu}(f_{\infty},f_{-\infty}) - 1) < 0 ,
\end{equation}
for some $k \ge 1$ and $\nu$,
then \eqref{nodal.eq} has a solution
$u \in T_{k,\nu}$.
\end{thm}

\begin{proof}
The proof is essentially the same as the proof of
\cite[Theorem~7.1]{RYN4}
(although the half-eigenvalues used in \cite{RYN4} were for problems of
the form \eqref{alt_form_hevals.eq}),
so we merely sketch it here.

By Theorem~\ref{branches.thm} there is a continuum, $\C_{k,\nu}$,
of solutions of \eqref{bif.eq},
bifurcating from $(\la_k/f_0,0)$,
and by the argument in \cite{RYN4} it can be shown that $\C_{k,\nu}$
`meets $(\la_{k,\nu}(f_{\infty},f_{-\infty}),\infty)$'
(more precisely, there exists a sequence
$(\mu_n,u_n) \in \C_{k,\nu}$, $n = 1,2,\dots,$
such that $\mu_n \to \la_{k,\nu}(f_{\infty},f_{-\infty})$,
$|u_n|_0 \to \infty$).
Since $\C_{k,\nu}$ is connected, it now follows from
\eqref{jumping_cond.eq}
that $\C_{k,\nu}$ must intersect the hyperplane
$\{1\} \X T_{k,\nu}$ at a point
$(1,u)$, and hence $u \in T_{k,\nu}$ is a solution of
\eqref{nodal.eq}.
\end{proof}

\begin{remarks}
Nodal solutions for similar problems have been obtained previously,
in several papers:
\begin{itemize}
\item
\cite{DR}, \cite{RYN3} and \cite{RYN4}
considered one separated and one multi-point Dirichlet boundary
condition;
\item
\cite{RYN5}
considered two multi-point Dirichlet boundary conditions;
\item
\cite{RYN6}
considered two multi-point Neumann boundary conditions, or a
mixture of Dirichlet and Neumann conditions.
\end{itemize}
NB the papers \cite{DR}, \cite{RYN5}, \cite{RYN6} dealt with
equations involving the $p$-Laplacian.
\\
We briefly summarize the result obtained in these papers.
\\[1 ex](a)\ \
Theorem~5.7 in \cite{DR} obtained similar nodal solutions to
Theorem~\ref{jumping.thm} above,
but the limits \eqref{f_limits.eq} were not assumed to exist
(and the half-eigenvalues had not been obtained) in
\cite{DR},
so instead of the above condition \eqref{jumping_cond.eq}
a condition involving certain $\limsup$'s and $\liminf$'s of $f(s)/s$,
as $s \to \infty$, was used.
When applicable, condition \eqref{jumping_cond.eq} is a weaker
hypothesis than the conditions in \cite[Theorem~5.7]{DR}.
\\[1 ex](b)\ \
Theorem~7.1 in \cite{RYN4} is similar to Theorem~\ref{jumping.thm}
above (using half-eigenvalues), but only considers a 3-point boundary
condition, at one end point.
\\[1 ex](c)\ \
When $f_{\infty} = f_{-\infty}$
the nonlinearity $f$  is `asymptotically linear'.
Theorem~\ref{jumping.thm} clearly applies to such problems,
and the half-eigenvalues
$\la_{k,\nu}(f_{\infty},f_{-\infty})$
reduce to $\la_k/f_\infty$.
Such problems were considered in all the above cited papers.
\\[1 ex](d)\  \
The `superlinear' case $f_{\infty} = f_{-\infty} = \infty$
has also been considered,
see \cite[Theorem~5.5]{DR} and \cite[Theorem~5.4]{RYN3}.
This case does not involve half-eigenvalues, and the proofs of these
results extend readily to the present setting,
so we simply state the corresponding result here, without proof.
\end{remarks}

\begin{thm} \label{superlinear.thm}
Suppose that  $f_{\infty} = f_{-\infty} = \infty$.
If  \eqref{sfs_pos.eq} holds and $\la_k/f_0  > 1$, for some $k_0 \ge 1$,
then for each $k \ge k_0$, \eqref{nodal.eq} has a solution
$u \in T_{k,\pm}$.
\end{thm}

\end{document}